\documentclass[a4paper]{article}
\usepackage[T1]{fontenc}
\usepackage[utf8]{inputenc}
\usepackage{geometry}
\usepackage{color}
\usepackage{graphicx}
\usepackage{pict2e}
\usepackage{amsmath, amsfonts, amssymb, amsthm}
\usepackage{mathrsfs}
\usepackage{url}
\usepackage{siunitx}
\usepackage{mathdefs}
\usepackage[margin=5pt,font=small,labelfont=bf]{caption}

\renewcommand{\transpose}{^\mathrm{T}}

\begin{document}
\title{A Numerical Slow Manifold Approach to Model Reduction for Optimal Control of Multiple Time
  Scale ODE} \author{Dirk Lebiedz\thanks{Institute for Numerical Mathematics, Ulm University,
    Helmholtzstr. 20, 89081 Ulm, Germany} and Marcel Rehberg\thanks{Institute for Numerical
    Mathematics, Ulm University, Helmholtzstr. 20, 89081 Ulm, Germany}} \date{\today} \maketitle

\begin{abstract}
Time scale separation is a natural property of many control systems that can be exploited,
theoretically and numerically. We present a numerical scheme to solve optimal control problems with
considerable time scale separation that is based on a model reduction approach that does not need
the system to be explicitly stated in singularly perturbed form. We present examples that highlight
the advantages and disadvantages of the method.
\end{abstract}

\section{Introduction}
Optimization based control in practice depends on accurate models with small prediction error and
computation of a (feedback) control that is close or at least consistent with the true optimal
control for the process under consideration \cite{Marquardt2002}. Often the desired accuracy can
only be provided by nonlinear large scale models which leads to problems in online control, for
example via nonlinear model predictive control (NMPC), \cite{Findeisen2003, Diehl2002b} due to the
computational demand. Model reduction therefore plays an important part in the development of
control systems, see the paper by Marquardt \cite{Marquardt2002} who gives a concise review
of model development and model reduction techniques.

Model reduction can be divided into model order reduction which aims at decreasing the dimension of
the state space and model simplification which tries to simplify the evaluation of the model
equations. Both approaches can be combined and essentially strive to capture the most important
features of the dynamic process at the cost of an error in the reduced compared to the full
model. The trade off between lost accuracy and benefits of the reduced model always has to be
considered and carefully balanced depending on the application at hand.

\section{Model Order Reduction}
In the remainder of this article we will only refer to model order reduction. Therefore we introduce
the system
\begin{equation} \label{eq:baseSys}
  \dot{\tilde{z}} = \tilde{f}(\tilde{z},u), \, \tilde{z}(0)=\tilde{z}_0
\end{equation}
with state $\tilde{z}(t) \in \R^n$ and control $u(t) \in \R^m$. The right hand side $\tilde{f}: \R^n
\times \R^m \rightarrow \R^n$ is assumed to be in $C^\infty$. A general approach to model order
reduction can be summarized in the following steps \cite{Marquardt2002}:
\begin{enumerate}
\item Find a diffeomorphism $T: \R^n \rightarrow \R^n$ that maps $\tilde{z}$ via
\[
\tilde{z}-\tilde{z}^* = T(z) \, \Leftrightarrow \, z = \tilde{z}^* + T^{-1}(\tilde{z}),
\]
onto the new state $z(t) \in \R^n$, where $\tilde{z}^*$ is a possibly nonzero set point. The aim
of this coordinate change is to separate directions in the phase space of \eqref{eq:baseSys} that
have strong contributions to the dynamics from those that only contribute in a minor way.
\item Decompose the new state space into $x(t) \in \R^p$ and $y(t) \in \R^q$ such that $z=(x,
  y)\transpose$ and $n=p+q$. Here $x$ will play the role of the dominant states.
\item Assemble new dynamic systems for $x$ and $y$ from
\[
\dot{z} = (\diff{z} T(z))^{-1} \tilde{f}(\tilde{z}^* + T(z),u)
\]
and obtain
\begin{spl}
  \dot{x} &= f(x,y,u), \quad x(0) = x_0, \\
  \dot{y} &= g(x,y,u), \quad y(0) = y_0.
\end{spl}
The smoothness of the right hand sides $f: \R^p \times \R^q \times \R^m \rightarrow \R^p$ and $g:
\R^p \times \R^q \times \R^m \rightarrow \R^p$ is determined by the smoothness of $T$ and $T^{-1}$. 
\item Eliminate the dynamic equation for $y$ by one of the following methods: 
  \begin{description} 
  \item[Truncation:] Set $y = 0$ for the reduced dynamic system
    \[
    \dot{\tilde{x}} = f(\tilde{x},0,u), \quad \tilde{x}(0) = x_0
    \]
    with $\tilde{x} \approx x$ and state space dimension $m$. 
  \item[Residualization:] Set $\dot{y} = 0$ to obtain the differential-algebraic system
    \begin{spl}
      \dot{\tilde{x}} &= f(\tilde{x},\tilde{y},u), \quad \tilde{x}(0) = x_0, \\
      0 &= g(\tilde{x},\tilde{y},u).
    \end{spl}
    The dimension of the model is not reduced.
  \item[Slaving:] Obtain a map $\tilde{y} = \phi(x,u)$ either from the residualization approach by solving
    the algebraic equation explicitly or through an independent method. Using
    \[
    \dot{\tilde{x}} = f(\tilde{x},\phi(\tilde{x},u),u), \quad \tilde{x}(0) = x_0
    \]
    leads to a reduced model with state space dimension $m$. 
  \end{description}
\end{enumerate}

Several model order reduction methods have been proposed in the past and we proceed to give a short
description of some of them. 

Nonlinear balancing is an analytical method based on the theory of nonlinear Hankel operators and
their attributed singular value functions aimed at obtaining a nonlinear map $T$
\cite{Fujimoto2010}. In practice, empirical balancing that incorporates samples of the systems
behavior for different inputs and initial values can be used \cite{Hahn2002}. In that case the map
$T$ is linear.

Proper orthogonal decomposition (POD) is based on sampling representative trajectories of
\eqref{eq:baseSys}, called snapshots \cite{Kerschen2005}. Similarly to balancing a linear
transformation matrix $T$ is obtained using singular value decomposition. For a focus in the context
of optimal control see \cite{Kunisch2008}.

Additional approaches include combinations of balancing and POD \cite{Lall2002} and moment matching
for nonlinear systems \cite{Astolfi2010}.

\subsection{Slow Invariant Manifolds}
Our approach to model order reduction is based on time scale separation which is a frequent feature
of complex dynamic processes. A theoretical environment for such system is provided in singular
perturbation theory \cite{Hoppensteadt1971} where we deal with a system of the form
\begin{spln} \label{eq:singPert}
  \dot{x} &= f(x,y,u, \varepsilon), \quad x(0) = \xi(\varepsilon), \\
  \varepsilon \dot{y} &= g(x,y,u,\varepsilon), \quad y(0) = \eta(\varepsilon).
\end{spln}
The parameter $\varepsilon$ is assumed to be small ($0 < \varepsilon \ll 1$) and reflects the time
scale separation. The fast modes $y$ evolve on the time scale $\bigo(\varepsilon^{-1})$ whereas the
slow state dynamics are $\bigo(1)$. We put forward the following assumptions:
\begin{enumerate}
\item[A1] The involved functions $f,g,\xi$, and $\eta$ are at least $R+2$ times continuously
  differentiable with respect to their arguments on their respective domains of interest.
\item[A2] Let $\varepsilon = 0$ in \eqref{eq:singPert}, then the reduced system is given by
  \begin{spl}
    \dot{x} &= f(x,y,u,0), \quad x(0) = \xi(0), \\
    0 &= g(x,y,u,0), \quad y(0) = \eta(0).
  \end{spl}
  There exist solutions $x(t)= x_0(t)$ and $y(t) = y_0(t)$ of the reduced system for $t \in [0,T]$.
\item[A3] The Jacobian
\[
g_y = \diff{y} g(x_0(t), y_0(t),u(t),0) \in \R^{q \times q}
\]
has $1 \leq k \leq q$ eigenvalues $\lambda_i,$ $i = 1,2,\dotsc,k$ with $\Re(\lambda_i) < -\mu$ and
$q-k$ eigenvalues with $\Re(\lambda_i) > \mu$, $i = k+1,k+2,\dotsc,q$ where $\mu > 0$. In other
words the Jacobian has no purely imaginary eigenvalues and there is at least one stable direction.
\end{enumerate}
From the third condition it follows that $g_y$ is nonsingular, so the algebraic equation $0 =
g(x,y,u,\varepsilon)$ can (at least locally) always be solved with respect to $y$.  Let the
assumptions A1--A3 hold, then there exists (\cite{Hoppensteadt1971}, Theorem 1) an $\varepsilon_0$
and a $k$-dimensional manifold $S(\varepsilon)$ such that the solution of \eqref{eq:singPert} can be
expanded into series representations for $\varepsilon < \varepsilon_0$ provided that
$\eta(\varepsilon) \subset S(\varepsilon)$:
\begin{spl}
    x(t,u, \varepsilon) &= x^*(t,u, \varepsilon) + X(t/\varepsilon,u, \varepsilon), \\
    y(t,u, \varepsilon) &= y^*(t,u, \varepsilon) + Y(t/\varepsilon,u, \varepsilon)
\end{spl} 
with
\begin{spl}
  x^*(t,u, \varepsilon) &= \sum_{r=0}^R x_r^*(t, u) \varepsilon^r + \bigo(\varepsilon^{R+1}),  \quad  
  y^*(t,u, \varepsilon) = \sum_{r=0}^R y_r^*(t, u) \varepsilon^r + \bigo(\varepsilon^{R+1}), \\
  X(t/\varepsilon, u, \varepsilon) &= \sum_{r=0}^R X_r(t/\varepsilon, u) \varepsilon^r  + 
  \bigo(\varepsilon^{R+1}), \quad 
  Y(t/\varepsilon, u, \varepsilon) = \sum_{r=0}^R Y_r(t/\varepsilon, u) \varepsilon^r + 
  \bigo(\varepsilon^{R+1}).
\end{spl}
The fast motions are captured in the so called boundary layer corrections $X(t/\varepsilon,u)$ and
$Y(t/\varepsilon,u)$ which converge to $0$ exponentially fast. For the purpose of model order
reduction we neglect the boundary layer correction and focus on the slow or outer solution
$x^*(t,u,\varepsilon)$. A central result is the following that goes back to Fenichel
\cite{Fenichel1979}, see also \cite{Kaper2002}, and is related to the geometric singular perturbation
approach to the problem.
\begin{theorem}[\cite{Kaper2002}, Theorem 2.1, Fenichel, asymptotically stable slow manifolds]
Let assumptions A1-A3 hold. Then, for any sufficiently small $\varepsilon$, there is a function $h$
that is defined on a compact domain $K \subset \R^q \times \R^m$ such that the graph
\[
\mathcal{M}_\varepsilon = \{(x,y) \setsep y=h(x,u,\varepsilon),\, (x,u) \in K \}
\]
is locally invariant under \eqref{eq:singPert}. The function $h$ admits an asymptotic
expansion,
\begin{equation} \label{eq:hExpand}
h(x,u, \varepsilon) = \sum_{r=0}^R  h_r(x,u)\varepsilon^r  + \bigo(\varepsilon^{R+1}).
\end{equation} 
\end{theorem}
The manifold $\mathcal{M}_\varepsilon$ is also called slow invariant manifold (SIM). Utilizing 
$h(x,u,\varepsilon)$ we can reduce the system \eqref{eq:singPert} to
\[
\dot{x}^* = f(x^*, h(x^*,u,\varepsilon), u, \varepsilon), \quad x(0) = \xi(\varepsilon),
\]
in accordance to the slaving approach introduced earlier. In practice only an approximation to $h$
can be feasibly computed. Using $h_0$ corresponds to setting $\varepsilon=0$ in \eqref{eq:singPert}
and solving the algebraic system for $y$, the fast states are assumed to be relaxed immediately.

Using an explicit formula or the SIM for model order reduction presumes that the system is given in
singularly perturbed form and therefore the method skips step 1 of the algorithm discussed
earlier. For general systems \eqref{eq:baseSys} a nonlinear coordinate transformation $T$ would have
to be explicitly known to apply the outlined theory. For systems where the small parameter
$\varepsilon$ can be identified, conditions on when such a transformation exists and how it might be
constructed are given in \cite{Marino1988}. Otherwise a state space decomposition into fast and slow
modes has to be based on physical insight or numerical methods. Among them are eigenvalue analysis
of the linearized system or singular value analysis of sensitivity matrices \cite{Lebiedz2005c}.

\subsection{Approximation of the SIM}
To approximate the SIM, in general only a numerical procedure will be feasible either because the
analytic computations to obtain the coefficients for the asymptotic expansion \eqref{eq:hExpand} of
$h$ can not be carried out explicitly or the system can not be transformed to the singular perturbed
form. In this section we therefore regard the general system
\begin{spl} 
  \dot{x} &= f(x,y,u),\\
  \dot{y} &= g(x,y,u).
\end{spl}
We assume that the state space decomposition into fast modes $x$ and slow modes $y$ is already
carried out either by employing a priori knowledge or by using one of the methods mentioned
above. Our approach \cite{Lebiedz2011a,Lebiedz2004c,Reinhardt2008,Lebiedz2010} is based on
optimization of trajectory pieces or points in the state space where the slow variables are fixed at
a certain point in time $t^*$ and the according fast states are computed in dependence of the slow
variables, i.e. the SIM is parametrized by the slow states and can represented by a function smooth
function $y = h(x,u)$. The underlying idea is that the fast states will relax onto the SIM as fast
as the system dynamics possibly allow and then stay on it. A parametrization for $u$ also has to be
chosen. We will use a piecewise constant control function later in a multiple shooting approach to
solve the optimal control problem, therefore we will use $u = \mathrm{const}$ here. The optimization
problem for the computation of a reduced model is
\begin{spln} \label{eq:redProb}
&  \min_{y} \Phi(x,y,u) \\
  \text{subject to: } & \dot{x} = f(x,y,u) \\
  & \dot{y} = g(x,y,u) \\
  & x(t^*) = x^*, \quad u = u^*. 
\end{spln}
Let $F$ be the full right hand side vector, hence
\[
F(x,y,u) = \begin{pmatrix} f(x,y,u) \\ g(x,y,u) \end{pmatrix}
\]
and $J$ be the Jacobian of the full dynamic system with respect to $x$ and $y$, i.e.
\[
J =  \diff{x,y}{} F(x,y,u)
\]
For $\Phi$ we will either use 
\begin{equation} \label{eq:moreInt}
\Phi(x,y,u) = \int_0^{t^*} \normeuk{J F}^2 \dt.
\end{equation}
or 
\begin{equation} \label{eq:moreLocal}
\Phi(x,y,u) = \normeuk{J F}^2.
\end{equation}
In the first objective the value of $x$ is fixed at the end of the integration interval because the
fast modes are unstable in backward time. This means that trajectories starting at points $y^* =
y(t^*)$ that are not on the manifold will exponentially move away from it and thus a large
contribution to the objective function is created. In the second case the dynamic equations are no
longer constraints of the optimization problem.
\begin{remark}
$\normeuk{J F}^2$ in both objective functionals can be linked to minimizing curvature in the phase
  space \cite{Lebiedz2010} and a variational principle \cite{Lebiedz2011a}. There is a relation to
  the zero derivative principle \cite{Zagaris2009}.
\end{remark}

Additionally, for the application of $h(x,u)$ in optimal control we also need at least first order
sensitivities $\diff{x} h(x,u)$ and $\diff{u} h(x,u)$. However, they can be easily obtained from the
KKT system at the solution point $y^*$ of the optimization problem \eqref{eq:redProb}
\cite{Siehr2012a}.

Numerically, for the integral based objective \eqref{eq:moreInt} either single shooting or
collocation is used to obtain a nonlinear program (NLP) which is subsequently solved with an
interior point method implemented in the software package \texttt{IPOPT} \cite{Waechter2006}. The
local formulation \eqref{eq:moreLocal} is solved with a general Gauß-Newton method specifically
tailored to the problem at hand \cite{Siehr2012a}. In all cases warm starts are employed to improve
convergence of the optimization if the problem has to be solved for a series of fixed values
$(x^*,u^*)$.

\section{Optimal Control and Reduced Models}
Singularly perturbed optimal control problems have been studied extensively in the past
\cite{Kokotovic1984, Vasileva1986, Naidu2002, Dmitriev2006}. The main idea is to decompose the full
system into a slow and a fast part and infer properties and solutions of the full problem from the
independent analysis of the two subproblems. The linear case is well understood, but for nonlinear
systems the situation is much more difficult and only for special cases useful explicit results can
be obtained. To highlight some of the problems and what we can expect at most from our approach we
review shortly the nonlinear state regulator problem \cite{OMalley1978}:
\begin{spln} \label{eq:singOpt}
  &  \min_{u} j(u) = E(x(1),\varepsilon y(1),\varepsilon) + \int_0^1 f_0(x,y,u,\varepsilon) \dt \\
  \text{subject to: } & \dot{x} = f(x,y,u,\varepsilon), \quad x(0) = x_0, \\
  & \varepsilon \dot{y} = g(x,y,u,\varepsilon), \quad y(0) = y_0, \\
  & x(t) \in \R^p, \quad y(t) \in \R^q, \quad u(t) \in \R^m.
\end{spln}
For convenience all functions are supposed to be $C^\infty$ functions of their arguments on any
domain of interest. The $\varepsilon y(1)$ in the final time cost avoids $E$ to depend on the fast
variable $y$ for $\varepsilon = 0$. Using the Pontryagin minimum principle with the Hamiltonian
\[
\mathcal{H}(x,y,\lambda_x,\lambda_y,u,\varepsilon) = f_0 + \lambda_x\transpose f +
\lambda_y\transpose g
\]
we get (additionally to the primal dynamic equations) the following ODE system for the adjoint
variables $\lambda_x$ and $\lambda_y$:
\begin{spln} \label{eq:singBVP}
\dot{\lambda}_x &= -\diff[]{x}\mathcal{H}, \quad 
\lambda_x(1) = \diff[]{x} E(x(1),\varepsilon y (1),\varepsilon), \\
\varepsilon \dot{\lambda}_y &= -\diff[]{y}\mathcal{H}, \quad 
\lambda_x(1) = \varepsilon \diff[]{y} E(x(1),\varepsilon y(1),\varepsilon).
\end{spln}
We note that there are no restrictions on the value of the control $u$ thus $\diff{u}\mathcal{H} =
0$ is a necessary condition for a minimum to occur. Moreover we assume the strong Legendre-Clebsch
condition, $\diff{uu}\mathcal{H}$ is positive definite, to hold. In that case a (locally) optimal
control $u(t)$ that minimizes the cost functional $j(u)$ exists \cite{Bryson1975} for $\varepsilon >
0$. It also allows to solve $\diff{u}\mathcal{H} = 0$ (locally) for
$u=\omega(x,y,\lambda_x,\lambda_y,\varepsilon)$ and replace it in \eqref{eq:singOpt} and
\eqref{eq:singBVP}. We now have a singularly perturbed boundary value problem. The main challenge
with problems of this type is to determine a reasonable reduced problem, i.e. setting
$\varepsilon=0$ in \eqref{eq:singOpt} and \eqref{eq:singBVP}. In general not all boundary values can
be satisfied and some of them have to be relaxed. In this case the choice is obvious: $y(0)$ and
$\lambda_y(1),$ the boundary values associated with the fast modes can not be fulfilled because for
$\varepsilon=0$ their values are determined by algebraic equations. Even if a reduced problem can be
stated, for nonlinear problems it is not generally possible to postulate conditions for a solution
to exist. Besides the smoothness assumptions the following restrictions are necessary.
\begin{enumerate}
\item[A1'] The reduced problem 
\begin{spl}
\dot{x} &= f(x,y,\omega,0), \quad x(0) = x_0, \\
\dot{\lambda}_x &= -\diff[]{x}\mathcal{H}(x,y,\lambda_x,\lambda_y,\omega,0), \quad 
\lambda_x(1) = \diff[]{x} E(x(1),0,0), \\
0 &= g(x,y,\omega,0), \\
0 &= -\diff[]{y}\mathcal{H}(x,y,\lambda_x,\lambda_y,\omega,0), \quad  \\
\end{spl}
has a unique solution $x^0(t)$, $y^0(t)$, $\lambda_x^0(t)$, and $\lambda_y^0(t)$.
\item[A2'] The Jacobian
\[
\mathcal{H}_y = \begin{pmatrix} \diff{\lambda_y y} \mathcal{H} & \diff{\lambda_y \lambda_y} \mathcal{H} \\
 - \diff{y y} \mathcal{H} & -\diff{y \lambda_y} \mathcal{H} \end{pmatrix} \in \R^{2q \times 2q}
\]
evaluated along $x^0(t)$, $y^0(t)$, $\lambda_x^0(t)$, and $\lambda_y^0(t)$ has no purely imaginary
eigenvalues, moreover we require it to have exactly $q$ eigenvalues with positive and $q$
eigenvalues with negative real part.
\end{enumerate}
For the problem at hand it turns out that $\mathcal{H}_y$ is block diagonal with two identical
blocks with opposite sign and symmetric and therefore the second condition on $\mathcal{H}_y$ is
automatically fulfilled if all eigenvalues have nonzero real part. The second condition guarantees
the solvability of the algebraic part of the reduced system and the stability of the boundary layer
corrections. If A1' and A2' hold, the solution of the full problem converges to the solution of the
reduced problem for $\varepsilon \rightarrow 0$ and the following series representations can be
stated \cite{OMalley1978, Hoppensteadt1971}:
\begin{spln} \label{eq:BVPSol}
  x(t,\varepsilon) &= x^*(t,\varepsilon) + X_L(t/\varepsilon,\varepsilon) +
  X_R(s/\varepsilon,\varepsilon), \\
  y(t,\varepsilon) &= y^*(t,\varepsilon) + Y_L(t/\varepsilon,\varepsilon) +
  Y_R(s/\varepsilon,\varepsilon), \\
  u(t,\varepsilon) &= u^*(t,\varepsilon) + U_L(t/\varepsilon,\varepsilon) + 
  U_R(s/\varepsilon,\varepsilon), \\
  j(u) &= j^*(t,\varepsilon) + J_L(t/\varepsilon,\varepsilon) + 
  J_R(s/\varepsilon,\varepsilon),
\end{spln}
with $s = 1-t$. Boundary layer corrections emerge at both ends of the time interval and appropriate
series representations can be found to all right-hand-side terms in \eqref{eq:BVPSol}. The
eigenvalue condition on $\mathcal{H}_y$ is necessary for the stability of the left and right hand boundary
layer corrections to be stable in forward backward time, respectively.

The manifold $h(x,u,\varepsilon)$ \eqref{eq:hExpand} is an intrinsic property of a singularly
perturbed system and it exists independently from the use of the system as constraint of an optimal
control problem. Hence it can be used to reduce the dimension of the optimal control problem
\eqref{eq:singOpt} by replacing $y$ and we find
\begin{spl}
  &  \min_{u} j(u) = E(x(1),\varepsilon h(x(1),u(1),\varepsilon),\varepsilon) + 
  \int_0^1 f_0(x,h,u,\varepsilon) \dt \\
  \text{subject to: } & \dot{x} = f(x,h,u,\varepsilon), \quad x(0) = x_0, \\
  & x(t) \in \R^p, \quad u(t) \in \R^m.
\end{spl}
Its solution will correspond to $x^*(t,\varepsilon)$, $u^*(t,\varepsilon)$, and $j^*(\varepsilon)$
which means we have no way of obtaining information about the boundary layer corrections in this
case. If the manifold is only an approximation of order $k$ with respect to its $\varepsilon$ series
representation the solutions of the reduced problem will also be approximations of order $k$.

\section{Numerical Solution of the Optimal Control Problem}
We will provide a short overview of the numerical methods we use to solve general optimal
control problems of the type
\begin{spl}
  &  \min_{x,u,T,p} j(x,u,p) \\
  \text{subject to: } & \dot{x} = f(x,u,r) \\
  & e(t,x,u,T,r) = 0 \\
  & i(t,x,u,T,r) \leq 0
\end{spl}
where $x(t) \in \R^m$ is the state, $u(t) \in \R^r$ is the control, $T \in \R$ is the final time,
and $r \in \R^s$ are parameters. We do not discuss the solvability of the problem and conveniently
assume that local solutions exist. In order to solve the problem numerically we have to discretize
the control and state functions $u(t)$ and $x(t)$, respectively. For this purpose we use multiple
shooting \cite{Bock1984}, which means we divide the overall time interval $[0,T]$ into $N$
subintervals with node points $t_k$, $k=0,1,\dotsc,N$, $t_k<t_{k+1}$, $t_0 = 0$, $t_N = T$. On each
interval the control is kept constant with values $u_k \in \R^m$. This constant input is used to
solve the ODE for $x(t)$ on each interval with the help of a numerical integration
routine. Introducing $x_k^0$, $k = 1,2,\dotsc,N$ as initial values for the dynamic equation on each
interval with states $x_k(t)$, $t \in [t_{k-1}, t_k]$ and $u_k$ as value of the control function we
can formulate a finite dimensional nonlinear program
\begin{spl}
  &  \min_{x_k^0,u_k,T,p} \sum_{k=1}^N j(x_k,u_k,p) \\
  \text{subject to: } & \dot{x_k} = f(x_k,u_k,p), \quad k = 1,2,\dotsc,N, \\
  & x_k(t_k) - x_{k+1}^0 = 0, \quad k = 1,2,\dotsc,N-1, \\
  & e(t,x_k,u_k,T,p) = 0, \quad t \in [t_{k-1},k], \; k= 1,2,\dotsc,N, \\
  & i(t,x_k,u_k,T,p) \leq 0 \quad t \in [t_{k-1},k],\; k= 1,2,\dotsc,N. \\
\end{spl}
The equality constraints $x_k(t_k) - x_{k+1}^0$ ensure the continuity of the solution at the
multiple shooting nodes. The method is implemented as a C++ program using \texttt{IPOPT}
\cite{Waechter2006} for solving the NLP, CppAD, a tool for automatic differentiation \cite{Bell2010}
for obtaining accurate derivatives in the NLP as well as for a BDF integrator \cite{Skanda2012} that
is used to solve the ODEs on the multiple shooting intervals and provide sensitivities.

\section{Evaluation of the Manifold}
Eventually, we have to evaluate the manifold map $y = h(x,u)$ for arbitrary points $(x,u) \subset
\R^p \times \R^m$. We do not regard $\varepsilon$ as argument of $h$ here, since it is either a
fixed parameter for the numerical solution of the optimal control problem in case of singularly
perturbed problems or we regard general problems without explicit dependence on a small parameter
$\varepsilon$. We will discuss two alternatives: Solving the model reduction problem online,
i.e. whenever an evaluation of $h(x,u)$ is needed while solving the optimal control problem, and
interpolation of offline precomputed data obtained by evaluating $h(x,u)$ on a discrete set of
points $\mathcal{C} \subset \R^p \times \R^m$. Both methods have intrinsic advantages and
disadvantages. The online method can be easily applied since no preparation steps have to be
taken. However, calculating $h(x,u)$ is costly and involves the solution of a NLP which could slow
down the overall computation. On the contrary the interpolation approach provides a direct, fast,
and easy to evaluate object, that promises a larger speed up. However, it suffers from the need to
precompute the manifold data and building the interpolation object, both tasks take a considerable
amount of time. This makes this approach only effective if the optimal control problem has to be
solved very often (e.g. in NMPC) so that the time spend in the preliminary stages is outweighed by
the overall performance gain. Furthermore, especially for higher dimensional problems, the storage
needed for the interpolation data might be to large to be handled properly for given hardware
resources.
\subsection{Online Evaluation}
If an evaluation of $h(x,u)$ is needed for a certain $(x_0,u_0)$ we solve the model reduction
problem \eqref{eq:redProb} with the slow states and control fixed to $(x_0,u_0)$. For performance
reasons only the local formulation \eqref{eq:moreLocal} is feasible because integration of the full
model is dispensed with and a general Gauß-Newton method can be used for efficient solution of the
minimization problem \cite{Siehr2012,Siehr2012a,Lebiedz2012}. Although it seems that this approach
contradicts the purpose of model reduction, since the full right hand side still has to be
evaluated, there is a computational advantage due to the decreased stiffness of the reduced
model. In addition the points at which $h(x,u)$ has to be evaluated will typically be close to each
other which makes it possible to use warm starts for the Gauß-Newton procedure. Thus, in general,
only very few iterations will be needed to solve the model reduction problem. This is in principle
similar to solving an explicitly given differential-algebraic equation where usually inside the
integration routine in each time step only a few iterations of a nonlinear equation solver are
needed for the algebraic part of the dynamic problem.

\subsection{Interpolation}
The interpolation approach is based on an interpolation function $\hat{h}(x,u)$ that will be used
instead of a pointwise computation of $h(x,u)$. We choose radial basis function (RBF) interpolation
\cite{Wendland2005} because it is independent of the dimension of the input space and configuration
of the interpolation nodes (i.e. grid free). The following short presentation of the subject is
strongly based on \cite{Wendland2005}. Given a set of nodes $\mathcal{C} = \{x_k\}_{k=1}^N$,
$\mathcal{C} \subset \Omega \subset \R^m$, radial basis interpolants are of the form
\[
s(x) = \sum_{k = 1}^N \lambda_k \Phi(x,x_k), \quad \lambda_k \in \R
\]
with basis functions $\Phi: \R^m \times \R^m \rightarrow \R$. The domain $\Omega \subset \R^m$ is
assumed to be open and bounded and satisfy an interior cone condition. The interpolation is carried
out by determining the coefficients $\lambda_k$ such that
\begin{equation}  \label{eq:genInterp}
L_k(s) =  L_k(f) = f_k, \quad f_k \in \R,
\end{equation}
holds, where $f,s \in H$ are functions, $H$ is a Hilbert space of functions $\R^d \rightarrow
\R^l$ and $L_k \in H^*$ are linear functionals from the dual of $H$. The following
property of a function $\Phi$ is central.
\begin{definition}[Positive definite function]
A continuous function $\Phi: \R^m \times \R^m \rightarrow \R$ is said to be positive definite if for
all $N \in \N$, all sets of pairwise distinct nodes $\mathcal{C} = \{x_k\}_{k=1}^{N}$, and all
$\alpha \in \R^N \setdiff \{0\}$ it holds that
\[
\sum_{\ell = 1}^N \sum_{k = 1}^N \alpha_\ell\alpha_k \Phi(x_\ell,x_k) > 0.
\]
\end{definition}
\begin{remark}
  Positive definite functions $\Phi(x,y)$ are also known as kernels and give rise to reproducing
  kernel Hilbert spaces, a topic we do not want to pursue any further here, see the book by Wendland
  \cite{Wendland2005} for more details.
\end{remark}
In the most simple case of pointwise evaluation as the functionals in \eqref{eq:genInterp},
i.e. $L_k (s) = \delta_{x_k} s = s(x_k)$, the interpolation condition becomes
\[
\delta_{x_\ell} s = s(x_\ell) = \sum_{k = 1}^N \lambda_k \Phi(x_\ell,x_k) = \delta_{x_\ell} f =
f_\ell
\]
and we obtain the linear system
\[
A \lambda = F, \quad A_{\ell,k} = \Phi(x_\ell,x_k), \quad 
\lambda = \begin{pmatrix} \lambda_1 &  \lambda_2 & \dotsc &\lambda_N \end{pmatrix}\transpose, \quad
F = \begin{pmatrix} f_1 & f_2 & \dotsc & f_N \end{pmatrix}\transpose,
\]
with $A$ positive definite since $\alpha\transpose A \alpha > 0$ holds for all $\alpha \in \R^N\setdiff
\{0\}$ by definition. That guaranties a unique solution to the interpolation problem for all sets of
pairwise distinct nodes.

The same approach can also used for Hermite interpolation where we interpolate not only the function
value itself but also partial derivatives. The relevant functionals $L_k$ are thus given by
\[
L_k = \delta_{x_k} \circ \diff[\alpha_k]{} , \quad k = 1,2,\dotsc,N, \quad \alpha_k \in \N^d,
\]
where the differential operator $\diff[\alpha_k]{}$ indicates concatenated derivatives according to
the multi-index $\alpha_k$. In general we demand that $x_k \neq x_\ell$ or $\alpha_k \neq
\alpha_\ell$ for $k \neq \ell$ to guarantee linear independence of the functionals and therefore a
unique solution to the interpolation problem. With our model reduction procedure for each
interpolation node the function value and all partial derivatives of first order can be computed. If
the interpolation nodes are chosen pairwise distinct than the linear independence follows. The
interpolant is given by
\[
s(x) = \sum_{k = 1}^N \lambda_k \diff[\alpha_k]{2} \Phi(x,x_k).
\]
The subscript $2$ of the differential operator indicates differentiation with respect to the second
variable. The interpolation matrix has entries of the form 
\[
\diff[\alpha_\ell]{1}\diff[\alpha_k]{2} \Phi(x_\ell,x_k)
\]
and it can be shown that for positive definite and sufficiently smooth $\Phi$ it is again also
positive definite and thus the interpolation problem can be uniquely solved.

In practice univariate radial basis functions are commonly used, i.e. $\Phi(x,y) \define
\phi(\normeuk{x-y})$. We will use the Gaussian function
\[
\phi(r) = \e^{-c^2r^2}, \quad c \in \R, c > 0.
\]
The parameter $c$ is called the shape parameter. It plays an essential role for the interpolation
error and the stability of the interpolation process by virtue of the fact that it is strongly
connected to the condition of the interpolation matrix. In the case of pointwise interpolation the
entries of $A$ are $\e^{-c^2 \normeuk{x_\ell - x_k}}$, $\ell,k = 1,2,\dotsc,N$ which, for $c
\rightarrow 0$, will converge to $1$ for all $x_\ell, x_k \in \R^m$. Conversely, for $c \rightarrow
\infty$ the basis function $\phi$ will either converge to $0$ if $x_\ell \neq x_k$ or to $1$ if
$x_\ell = x_k$. This means the matrix $A$ tends to being singular in the first case and becoming the
unit matrix in the second case. The interpolation error $\norm{s-f}$ is also subject to the same
trade of principle. Determining a ``good'' $c$ is crucial for the performance of the
interpolation. Therefore we use the algorithm suggested in \cite{Rippa1999} which is based on a
computational favorable reformulation of a leave-one-out optimization scheme.

Since we want to use the interpolation in the optimization we are interested in a fast
evaluation. We use a partition of unity approach to divide the domain of interest $\Omega$ into
smaller subdomains and thereby bound the computing cost. To be more precise we look for a
overlapping covering of $\Omega$ by open and bounded sets $\Omega_j$, $j=1,2,\dotsc,M$ and
continuous functions $\omega_j(x): \R^m \rightarrow \R$ such that
\[
\sum_{j=1}^M \omega_j(x) = 1, \; \forall x \in \Omega \text{ and } \omega_j(x) = 0, \; \forall x \notin \Omega_j.
\]
Assuming that we have a feasible covering $\{\Omega_j\}$ we can build local interpolants $s_j(x)$
for each $\Omega_j$. Additionally let 
\[
I(x) = \{j | x \in \Omega_j\}
\]
be an index function that returns the indices of the patches a point $x$ is contained in. A global
interpolant is then simply given by
\[
s(x) = \sum_{j \in I(x)} \omega_j(x) s_j(x).
\]
Under certain conditions on the covering and the $\omega_j$ it can be shown that the global
interpolant enjoys the same error rates as in the naive global approach.

Two key ingredients are needed to really take computational advantage of the method: Partition
$\Omega$ in a way so that all $\Omega_j$ contain about the same amount of node points and have the
index function $I(x)$ be $\bigo(1)$ in terms of computing time, i.e. the cost of finding the patches
a random point lies in does neither depend on the overall number of centers $N$ nor on the number of
patches $M$. In that case evaluation is $\bigo(1)$ also, since for any number of nodes we can
partition $\Omega$ in a way that the number of points in a patch is below a certain constant
threshold which means that the sums, that have to be evaluated for each patch have a constant number
of terms which together with the constant time look-up leads to constant evaluation time.

If we assume the node points $x_k$ to be uniformly distributed in $\Omega$ one practical way to
achieve both aims is to use a fixed-grid structure which consists of axis parallel overlapping
boxes. Because of the parallelism many operations, like index querying can be independently
performed in each dimension. Given a set of nodes $\mathcal{C}$, overlap factor $\gamma \in
(0,0.5)$, and a lower bound for the average number of points in one box the number of boxes and
their border coordinates can be computed. If $o_i$, $i=1,2,\dotsc,d$ is the length of a box in
dimension $i$ the overlap factor $\gamma$ determines the fraction of $o_i$ by which two boxes
overlap in dimension $i$, $\gamma < 0.5$ ensures that a point $x$ can only be in maximal $2$ boxes
per coordinate direction.

The last missing part are the $w_j$. We choose dimension independent radial polynomials as suggested
in \cite{Tobor2004} for the same purpose. They are given by
\[
\omega_j(x) = \begin{cases} p \circ b_j(x) & x \in \Omega_j, \\
  0 \quad & \text{else,}
\end{cases}
\]
with $p: \R \rightarrow \R$,
\[
p(r) = -6r^5 + 15 r^4 - 10r^3 + 1
\]
a polynomial that fulfills the spline like conditions $p(0) = 1$, $p(1) = 0$, $\diff[k]{}p(0) =
\diff[k]{}p(1) = 0$, $k = 1, 2$ and $b_j: \R^m \rightarrow [0,1]$,
\[
b_j(x) = 1 - \prod_{i =1 }^d \frac{4(x_i - l_i^j)(r_i^j - x_i)}{(r_i^j - l_i^j)^2},
\]
where $l^j \in \R^m$ and $r^j \in \R^m$ are the lower left and upper right corner coordinates of the
box $\Omega_j$, respectively. The polynomial $p$ gives rise to a two times continuously
differentiable global interpolant $s(x)$.

Besides function evaluation we at least also need first order partial derivatives during the
optimization procedure which are computed by differentiating the interpolation function
symbolically with respect to $x_i$.

\section{Results}
We present some numerical examples. Computation times given were obtained on an Intel Xeon E5620
(\num{2.40} GHz) machine running \num{64} bit Debian Squeeze.
\subsection{Enzym Kinetics}
The first example is based on a simple Michaelis-Menten enzyme kinetics in singularly perturbed form
\cite{Murray1993}. The (full) problem is
\begin{spln} \label{eq:enzymFull}
  & \min \int_0^5 -50y + u^2 \dt \\
  \text{subject to: } & \dot{x} = -x + (x + 0.5)y + u,\\
  & \varepsilon \dot{y} = x - (x + 1.0)y, \\
  & x(0) = 1, \quad y(0) = y_0.
\end{spln}
The control and the objective are artificial and not related to a realistic model scenario. The reduced
problem is obtained by replacing $y$ with $h(x,u)$, eliminating the ODE for $y$ and the
initial condition $y_0$. The reduced problem is thus
\begin{spln} \label{eq:enzymRed}
  & \min \int_0^5 -50 h(x,u) + u^2 \dt \\
  \text{subject to: } & \dot{x} = -x + (x + 0.5)h(x,u) + u,\\
  & x(0) = 1. 
\end{spln}
We set $\varepsilon=\num{e-2}$ and use the multiple shooting approach described above for numerical
solution. To obtain initial values for $x$ and $y$ at the multiple shooting nodes we integrate the
uncontrolled system $(u(t)=0)$ numerically starting at $x(0)=0$ and $y(0) = y_0$. The overall time
interval is divided into 40 equidistant multiple shooting intervals. The termination tolerance for
\texttt{IPOPT} was set to \num{e-4} and the integration tolerance of the BDF-integrator to
\num{e-6}. The discretized full problem has \num{204} variables whereas the reduced problem has
\num{163}. Lastly, bounds are introduced for the state variables and the control. We use the lower
bounds $x_l=y_l=u_l=0$ and the upper bounds $x_u=y_u=u_u=5.5$ on all multiple shooting nodes.

The performance of the full system \eqref{eq:enzymFull} depends on the initial value
$y_0$. For $y_0 = 0$ we have an average runtime of \SI{2.5}{\second} and \num{26} NLP
iterations. For $y_0 = 0.5$, which is the first order approximation $h_0(1,0)$ we find
\SI{2.1}{\second} for \num{22} NLP iterations and lastly for $y_0 = 1$ we get \SI{2.2}{\second} and
\num{22} iterations.

Next we used the online computation of $h$. The algorithm clocks in at \SI{2.1}{\second} and
\num{22} iterations, which means no performance gain compared to the full problem. However, the
number of NLP iterations in the subproblem of approximating the manifold is interesting. The maximum
is \num{4} iterations but \SI{60}{\percent} of the calls terminate after \num{2} and
\SI{37}{\percent} only after \num{1} iterations. The call to the model reduction routine is by now
done through an external library which produces a lot of overhead, for example in terms of right
hand side function evaluations. Tight integration into the BDF integration algorithm might lead to a
significant speed up.

For the interpolation approach one first needs to choose a reasonable set of nodes. Although the RBF
method is grid independent for convenience we used a Cartesian grid on $[-0.5, -0.5] \times [10,
  10]$. The performance of the interpolator depends of course on the number of points but also on
the number of points per patch and overlap in the partition of unity approach. To assess the
influence we scatter searched the region $\{20,30,35,40\} \times \{0.025, 0.05,0.1,0.15\} \times
\{5,10,15\}$ for points in each direction, overlap and points per patch respectively. Median runtime
was \SI{0.87}{\second} and median number of NLP iterations \num{24.5}. Moreover, the fastest
combination took \SI{0.61}{\second} and only \num{21} iterations compared to the slowest which
needed \SI{1.4}{\second} and \num{23} iterations. It is apparent that using the interpolation
approach in this case is beneficial from an performance point of view. In the best case it is nearly
\num{4} times faster than the full problem. 

The main reason for the speed up is not so much the reduction in the number of optimization
variables but mainly the reduced stiffness along the manifold $h$. In the integration routine larger
step sizes are possible which greatly reduces the computational effort. This especially pays off in
the multiple shooting approach since the initial values at the multiple shooting nodes are subject
to optimization and they might be set away from the SIM for the full system in each iteration of the
NLP solver leading to transient behavior of the fast trajectories on each interval and therefore
forces the integrator to use small steps. An overview of example integrator statistics is given in
table \ref{tab:enzym} as well as the result of the performance tests.

\begin{table} 
  \begin{center}
    \caption{Summary of various statistics concerning the solution of problem \eqref{eq:enzymFull} and 
      \eqref{eq:enzymRed}. The steps statistics refer to the integration and are the sums over all NLP
      iterations and multiple shooting intervals.\label{tab:enzym}}
      \begin{tabular}{lrrrrr}
        \hline
        \multicolumn{1}{c}{Problem}
        & \multicolumn{1}{c}{Time} 
        & \multicolumn{1}{c}{NLP iter} 
        & \multicolumn{1}{c}{Time per iter} 
        & \multicolumn{1}{c}{Steps} 
        & \multicolumn{1}{c}{Rej. steps} \\
        \hline
        \eqref{eq:enzymFull} & \SI{2.3}{\second} & \num{23} & \SI{0.1}{\second} 
        & \num{43165} & \num{7369} \\
        \eqref{eq:enzymRed} online  & \SI{2.1}{\second} & \num{22} & \SI{0.1}{\second} 
        & \num{5520} & \num{1281} \\
        \eqref{eq:enzymRed} offline (median) & \SI{0.9}{\second} & \num{24.5} & \SI{0.04}{\second} 
        & -  & - \\
        \eqref{eq:enzymRed} offline (best) & \SI{0.6}{\second} & \num{21} & \SI{0.03}{\second} 
        & \num{5194} & \num{1236}\\
        \hline
      \end{tabular}
  \end{center}
\end{table}

\begin{figure}
 \begin{center}
 \parbox[t]{0.49\textwidth}{\centering
   \includegraphics{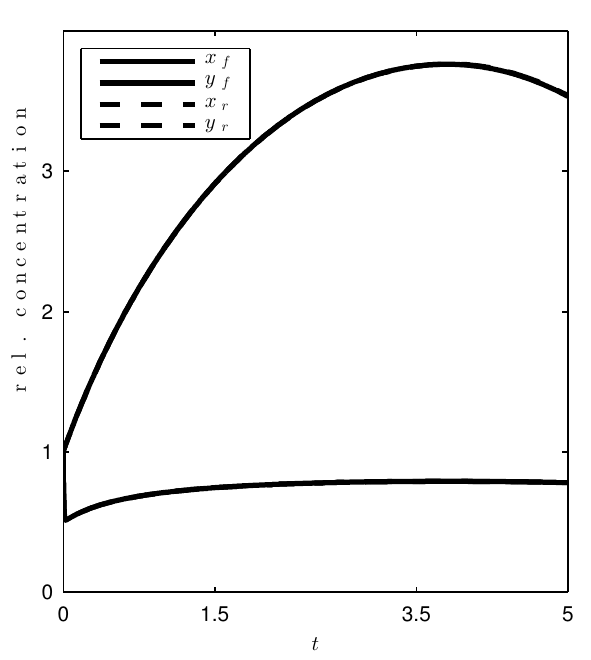}
   \caption{Example trajectories using the control from the full problem \eqref{eq:enzymFull},
     subscript $f$ and the reduced problem \eqref{eq:enzymRed}, subscript $r$. Both trajectories
     overlap.}
   \label{fig:enzymSpec}}
 \parbox[t]{0.49\textwidth}{\centering
   \includegraphics{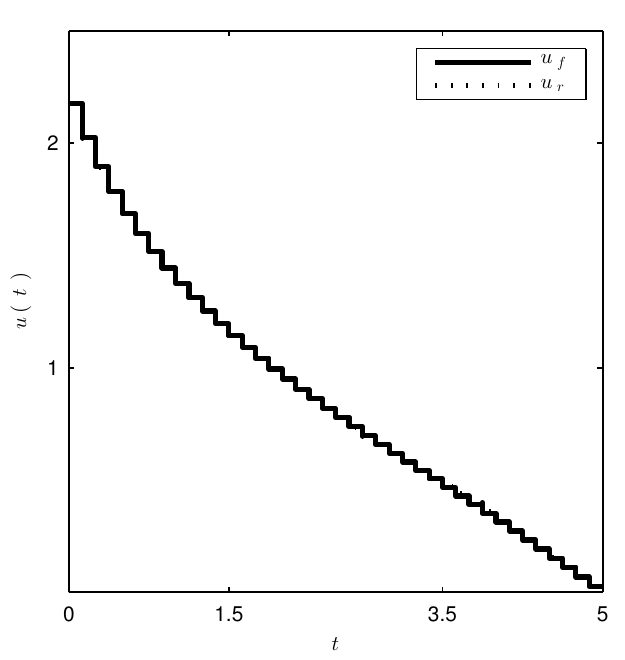}
   \caption{Example controls computed from the full problem \eqref{eq:enzymFull}, subscript $f$ and
     the reduced problem \eqref{eq:enzymRed}, subscript $r$. Both controls overlap.}
   \label{fig:enzymCon}}
 \end{center}
\end{figure}

If we use the computed optimal control from the reduced system as input of the full system we obtain
virtually the same objective values as for the control computed with the full system. In this
example the error of the reduction is negligible which can also be concluded from the system
itself. Example trajectories and controls are plotted in figure \ref{fig:enzymSpec} and
\ref{fig:enzymCon} respectively. In both plots the subscript $f$ refers to the results obtained
using the computed control from the full system \eqref{eq:enzymFull}, whereas $r$ refers to the
results using the control computed from the reduced problem \eqref{eq:enzymRed}.

\subsection{Voltage Regulator}
The next example is taken from \cite{Naidu1988}, example $4.2$. It describes a voltage regulator
with $5$ states governed by a set of linear ODEs. The problem is given by
\begin{spln} \label{eq:VRFull}
  & \min \frac{1}{2}\int_0^2 x_1^2 + u^2 \dt \\
  \text{subject to: } & \dot{x}_1 = -\frac{1}{5} x_1 + \frac{1}{2} x_2, \\
  & \dot{x}_2 = -\frac{1}{2} x_2 + \frac{8}{5} y_1, \\
  & \varepsilon \dot{y}_1 = - \frac{5}{7} y_1 + \frac{30}{7} y_2, \\
  & \varepsilon \dot{y}_2 = -\frac{5}{4} y_2 + \frac{15}{4} y_3, \\
  & \varepsilon \dot{y}_3 = -\frac{1}{2} y_3 + \frac{3}{2} u.
\end{spln}
The problem, although in essence linear-quadratic, has some interesting features: The coupling
between the slow and fast subsystems is only through the one fast state $y_1$ which means that only
one state has to be reproduced during the optimization, i.e. we are only interested in $y_1 =
h(x_1,x_2,u)$. The reduced problem is
\begin{spln} \label{eq:VRRed}
  & \min \frac{1}{2}\int_0^2 x_1^2 + u^2 \dt \\
  \text{subject to: } & \dot{x}_1 = -\frac{1}{5} x_1 + \frac{1}{2} x_2, \\
  & \dot{x}_2 = -\frac{1}{2} x_2 + \frac{8}{5} h(x_1,x_2,u).
\end{spln}
First, we set $\varepsilon=\num{0.2}$ and solve both problems on \num{10} multiple shooting
intervals with the \texttt{IPOPT} tolerance set to \num{e-3}. Again a selection of initial values
for the full system was used. An overview is given in table \ref{tab:VRInit}. Note, that this
selection leads to a set of two initial values for the reduced system, namely $x_0 =
(-10,0)\transpose$ and $x_0 = (-10,10)\transpose$. The initial values for the state variables at the
multiple shooting nodes are obtained through integrating the ode system with the initial control
$u(t)=0$. Bounds are introduced as follows: $y_1 \in [-20,20]$, $y_2 \in [10,50]$, $y_3, y_4, y_5 \in
[\num{-e8}, \num{e8}]$,  and $u \in [-15,15]$.

\begin{figure}
 \begin{center}
 \parbox[t]{0.49\textwidth}{\centering
   \includegraphics{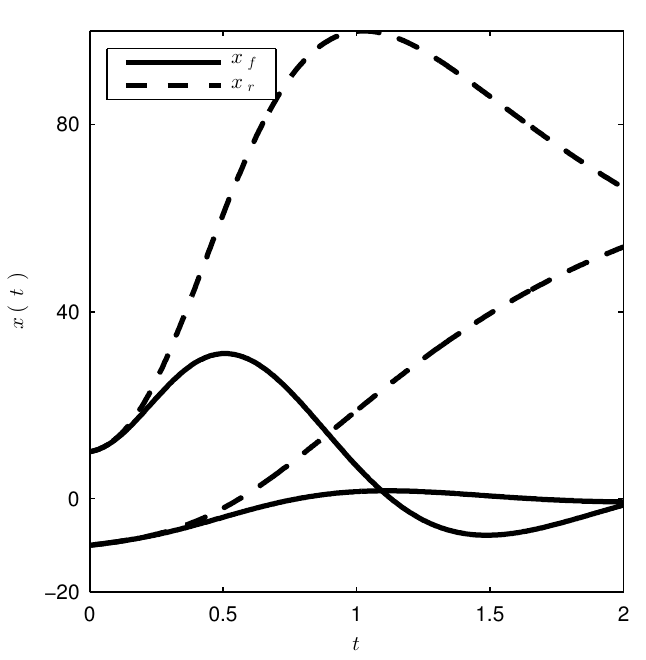}
   \caption{Example trajectories using the control from the full problem \eqref{eq:VRFull},
     subscript $f$ and the reduced problem \eqref{eq:VRRed}, subscript $r$ with $\varepsilon =
     0.2$.}
   \label{fig:VRSpec}}
 \parbox[t]{0.49\textwidth}{\centering
   \includegraphics{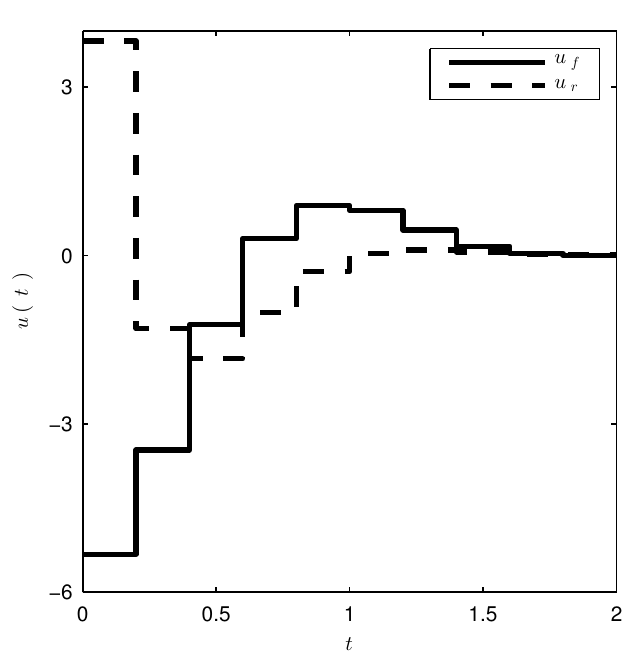}
   \caption{Example controls computed from the full problem \eqref{eq:VRFull}, subscript $f$ and
     the reduced problem \eqref{eq:VRRed}, subscript $r$ with $\varepsilon=0.2$.}
   \label{fig:VRCon}}
 \end{center}
\end{figure}

Using the control computed from the reduced problem for the full problem leads
to extremely large objective values in comparison and thus renders the model reduction unusable in
this case. This can also be seen in figures \ref{fig:VRSpec} and \ref{fig:VRCon} which compares the
full system with the reduced system solved with the online evaluation of $h(x,u)$. Additionally
there is no runtime advantage: The full system needs on average \SI{1.1}{\second} compared to the
offline approach which needs \SI{1.4}{\second} which is also the median timing of the online method.
\begin{table} 
  \begin{center}
    \caption{Final objective values of problems \eqref{eq:VRFull} and \eqref{eq:VRRed} for a
      selection initial values and $\varepsilon=0.2$. \label{tab:VRInit}}
      \begin{tabular}{lrrr}
        \hline
        \multicolumn{1}{c}{$x_0$}
        & \multicolumn{1}{c}{\eqref{eq:VRFull}} 
        & \multicolumn{1}{c}{\eqref{eq:VRRed}, online} 
        & \multicolumn{1}{c}{\eqref{eq:VRRed}, offline} \\
        \hline
        $(-10, 0, 0, 0, 0)$ & \num{32.9} & \num{35.1} & \num{34.9} \\
        $(-10, 0, 10, 0, 10)$ & \num{25.1} & \num{521.7} & \num{612.7} \\
        $(-10, 0, 0, 10, 10)$ & \num{24.7} & \num{648.0} & \num{750.1} \\
        $(-10, 10, 10, 10, 10)$ & \num{20.4} & \num{782.6} & \num{830.3} \\
        $(-10, 10, 0, 0, 10)$ & \num{20.5} & \num{521.4} & \num{559.4} \\
        $(-10, 10, 10, 0, 10)$ & \num{20.0} & \num{579.1} & \num{619.5} \\
        \hline
      \end{tabular}
  \end{center}
\end{table}

\begin{figure}
 \begin{center}
 \parbox[t]{0.49\textwidth}{\centering
   \includegraphics{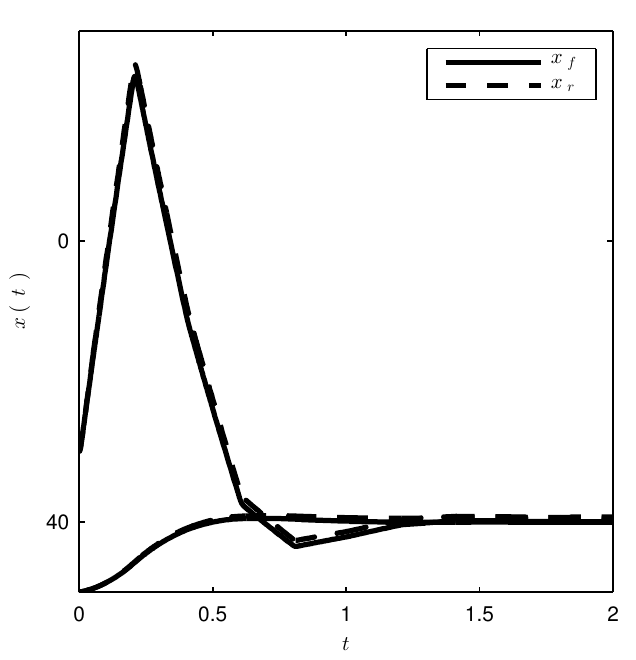}
   \caption{Example trajectories using the control from the full problem \eqref{eq:VRFull},
     subscript $f$ and the reduced problem \eqref{eq:VRRed}, subscript $r$ with $\varepsilon =
     \num{2e-3}$.}
   \label{fig:VRSpecEps2}}
 \parbox[t]{0.49\textwidth}{\centering
   \includegraphics{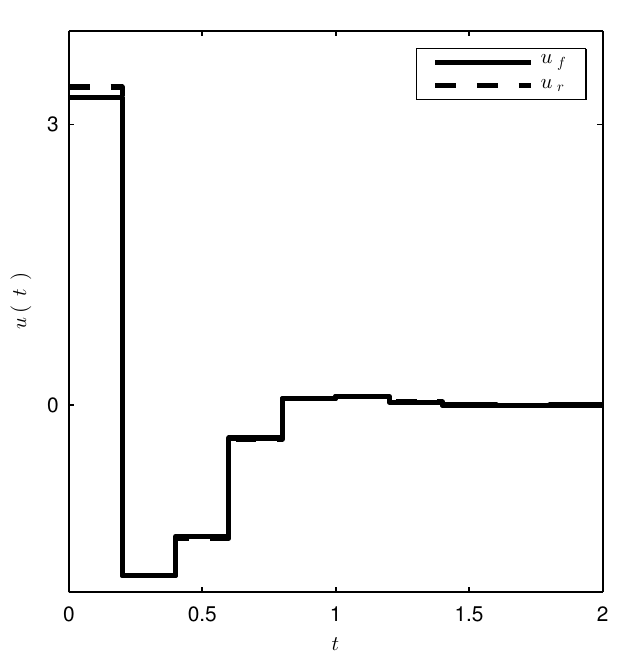}
   \caption{Example controls computed from the full problem \eqref{eq:VRFull}, subscript $f$ and
     the reduced problem \eqref{eq:VRRed}, subscript $r$ with $\varepsilon=\num{2e-3}$.}
   \label{fig:VRConEps2}}
 \end{center}
\end{figure}

If we increase the spectral gap by setting $\varepsilon = 2\cdot10^{-3}$ the results are much more
favorable. First of all the computed input from the reduced model is very close to solution of the
full problem. Therefore also the objective values are similar. See figures \ref{fig:VRSpecEps2} and
\ref{fig:VRConEps2} for an example. With the reduced model we find a significant
computational advantage, as documented in table \ref{tab:VRTimings}. The runtime for the full
problem depends strongly on the initial values and varies between \SI{2.3}{\second} to
\SI{5.1}{\second} which is between \num{5} to \num{10} times slower compared to the fastest solution
of the reduced problem with the offline method. The online approach is around \num{2} to \num{3}
times faster. A further advantage worth mentioning in this context is that the time needed for the
reduced problem is much less dependent on the initial values, which makes the computation more
reliable in online control scenarios where the next input has to be computed within a given time
frame.

As in the enzyme example we systematically tried various parameter combinations (points per
dimension, overlap and points per patch) for the interpolator. We already mentioned the best and
median runtime values, however it should also be noted that bad parameter combinations can decrease
the algorithmic performance significantly. The maximum time needed for both $\varepsilon$ values and
sets of initial values was over \SI{16}{\second}. In a considerable number of cases the problem
could not even be solved. This shows that the interpolator approach has to be tuned carefully but
further analysis reveals that at least in this case the best configuration is the same for both
initial values and $\varepsilon$.

\begin{table} 
  \begin{center}
    \caption{Summary of various statistics concerning the solution of problem \eqref{eq:VRFull} and
      \eqref{eq:VRRed} for $\varepsilon = \num{0.2 e-3}$. Timings are averages over all initial
      values. \label{tab:VRTimings}}
      \begin{tabular}{lrrr}
        \hline
        \multicolumn{1}{c}{Problem}
        & \multicolumn{1}{c}{Time} 
        & \multicolumn{1}{c}{NLP iter} 
        & \multicolumn{1}{c}{Time per iter} \\
        \hline 
        \eqref{eq:VRFull} & \SI{3.7}{\second} & \num{37.7} & \SI{0.1}{\second} \\
        \eqref{eq:enzymRed} online  & \SI{1.6}{\second} & \num{16.5} & \SI{0.1}{\second} \\
        \eqref{eq:VRRed} offline (median) & \SI{1.3}{\second} & \num{18} & \SI{0.07}{\second} \\
        \eqref{eq:VRRed} offline (best) & \SI{0.5}{\second} & \num{16} & \SI{0.03}{\second} \\
        \hline
      \end{tabular}
  \end{center}
\end{table}

\section{Concluding Remarks}
Given an optimal control problem, the aim of model reduction is to determine and solve a smaller
problem and use its solution as input to the full scale problem hoping for computational benefits
while still obtaining a feasible and nearly optimal control. Our approach to reduce the dimension of
the state space is based on time scale separation, i.e. processes evolving on slow and fast time
scales within the same system. According to \cite{Lebiedz2004c,
  Lebiedz2011a,Reinhardt2008,Siehr2012a} we formulate a nonlinear optimization problem that
identifies a slow manifold in the state space, parametrized by the slow states. This manifold hence
defines a map $y=h(x,u)$ of the slow variables and control onto the fast variables and can be used
to reduce the dynamic system.

Singular perturbation theory delivers a framework for optimal control problems involving fast/slow
differential systems. Using the Pontryagin minimum principle one arrives at a singularly perturbed
boundary value problem. Its solution consists of three components: A slowly varying part which
represents the system confined to the slow manifold and two fast vanishing boundary layer
corrections. Because we use an approximation to the slow manifold to represent the fast states and
thereby reduce the dimension of the state space we are only able to obtain an approximation to the
slow part of the optimal control solution. Two examples, both of which are singularly perturbed
systems, are used to illustrate that if the boundary layer corrections are small, the solution of
the reduced system can produce controls that drive the full system in a nearly optimal
fashion. Solving the reduced problem is up to ten times faster if the offline, interpolation based
method is used for the manifold map $h$. The numerical scheme presented here can be used unmodified
to solve general nonlinear optimal control problems, i.e. problems not explicitly in singular
perturbed form.

If larger systems, especially with more slow states are considered the interpolation approach will
suffer from the curse of dimensionality because of the exponentially growing number of nodes and
with it also interpolation data that has to be handled. This problem could to a certain degree be
overcome by using a more suited data structure (e.g. kd-trees) and by using a tight, problem
specific state space and control domain and take advantage of the ability of using scattered
nodes.

The online method is not subject to the dimensionality problem, however its evaluation for one input
point takes considerably longer since the full system has to be evaluated in the general Gauß-Newton
method. Still the benefit from reduced stiffness can speed up the overall solution of the optimal
control problem. A tight integration into the ODE integration routine, similar to a DAE solver would
greatly decrease unnecessary overhead and increase speed and stability.

\bibliography{/home/mrehberg/literature/literature} \bibliographystyle{plain}
\end{document}